\documentclass{amsart}

\usepackage{amsmath,amssymb,amsthm}
\usepackage[numbers]{natbib}

\newtheorem{theorem}{Theorem}[section]
\newtheorem{proposition}[theorem]{Proposition}
\newtheorem{lemma}[theorem]{Lemma}
\newtheorem{corollary}[theorem]{Corollary}
\theoremstyle{definition}
\newtheorem{definition}[theorem]{Definition}
\newtheorem{example}[theorem]{Example}
\theoremstyle{remark}
\newtheorem{remark}[theorem]{Remark}
\numberwithin{equation}{section}

\title[Rigidity at Isotropy Order $n-6$]{Rigidity of Constant-Curvature
Minimal Two-Spheres in Complex Hyperquadrics at Isotropy Order $n-6$}

\author{Zonghao Jie}
\address{School of Mathematical Sciences,
University of Chinese Academy of Sciences}
\email{jiezonghao23@mails.ucas.ac.cn}
\author{Xiaoxiang Jiao}
\address{School of Mathematical Sciences,
University of Chinese Academy of Sciences}
\email{xxjiao@ucas.ac.cn}

\begin{document}

\begin{abstract}
We classify linearly full, irreducible, totally unramified conformal minimal
two-spheres of constant Gaussian curvature in the real Grassmannian
$G(2,n;\mathbb R)$, or equivalently in the complex hyperquadric $Q_{n-2}$,
whose finite isotropy order is $n-6$, where $n\geq7$. We remove the
additional polynomiality assumption on a Frenet norm in the earlier
classification of Li and Jiao. An exterior-product norm identity reduces
the problem to the multiplicity of the irreducible polynomial
$1+zw$; a further divisibility argument handles the boundary dimension
$n=7$, including isotropy order one. The dimension $n$ is necessarily odd,
and the immersion is congruent, after a conformal reparametrization, to
the second member of a Veronese sequence with its prescribed real
structure. Its curvature is $2/(5n-13)$ in the Grassmannian normalization.
An explicit model in $Q_5$ illustrates the boundary case and its harmonic
sequence.
\end{abstract}

\maketitle

\section{Introduction}

The classification of minimal two-spheres of constant Gaussian curvature
is particularly rigid in complex projective space. The harmonic-sequence
classification of Bolton, Jensen, Rigoli and Woodward
\citep{bolton1988conformal} identifies the linearly full examples with
members of the Veronese sequence. In Grassmannians, additional geometric
hypotheses are needed. Harmonic sequences provide a framework in which
reducibility, isotropy order and ramification can be distinguished
\citep{burstall1986construction,chern1987harmonic}.

The complex hyperquadric $Q_{n-2}\subset\mathbb{CP}^{n-1}$ is the oriented
double cover of the real Grassmannian $G(2,n;\mathbb R)$. Harmonic
two-spheres in these spaces were studied by Wolfson
\citep{wolfson1986harmonic}, and Bahy-El-Dien and Wood
\citep{bahy1989explicit} constructed the finite-isotropy examples by
successive forward and backward replacements. The low-dimensional
classification program addresses the successive targets $Q_3$
\citep{li2016q3}, $Q_4$ \citep{jiao2021q4} and $Q_5$
\citep{jiao2019q5}, under the respective hypotheses of these works.
In general dimension, Jiao and Li \citep{jiao2017classification} treated
the irreducible, totally unramified case at the threshold $r\geq n-5$.
The next threshold, $r=n-6$, permits a two-step replacement construction.

Li and Jiao \citep[Proposition~4 and Theorem~1]{li2023conformal} obtained
the Veronese second-member classification at this threshold under the
additional assumption that the squared norm $|f_2^{(n-1)}|^2$ of the
associated Frenet vector is polynomial in $z,\bar z$, for the polynomial
lift of the initial holomorphic curve. Their Proposition~3 excludes the
one-step construction without this extra assumption. The purpose of
the present paper is to remove the polynomiality assumption in the
remaining two-step case.

We write $\underline v$ for the line spanned by a nonzero vector $v$.
Let $V_j^{(m)}$ denote the standard Veronese Frenet vectors, defined in
Section~\ref{sec:projective}, and put
\[
W_n=\operatorname{antidiag}(1,-1,\ldots,1).
\]
The metric and total-unramification conventions are specified in
Section~\ref{sec:preliminaries}. In particular, all curvatures below use
the Grassmannian metric, which pulls back to twice the standard
Fubini--Study induced metric on the quadric.

\begin{theorem}\label{thm:main}
Let $n\geq7$, and let $\phi:\mathbb S^2\to G(2,n;\mathbb R)$ be a
linearly full, irreducible, totally unramified conformal minimal
immersion of constant Gaussian curvature $K$ and finite isotropy order
$r=n-6$. Then $n$ is odd. After a conformal reparametrization of
$\mathbb S^2$ and an isometry of the target,
\begin{equation}\label{eq:main-model}
\phi=\overline{\underline{UV_2^{(n-1)}}}
       \oplus\underline{UV_2^{(n-1)}},\qquad
U\in U(n),\quad U^TU=W_n,\qquad K=\frac{2}{5n-13}.
\end{equation}
Conversely, for every odd $n\geq7$ the displayed models satisfy these
hypotheses, and all choices of $U$ satisfying $U^TU=W_n$ give the same
target-isometry class. No polynomiality assumption on
$|f_2^{(n-1)}|^2$ is required.
\end{theorem}

The new step is an identity in an exterior-power space. For a polynomial
lift $f_0$ of the initial holomorphic curve, put
$F_1=f_0\wedge f_0'$ and $F_2=f_0\wedge f_0'\wedge f_0''$.
The two-step construction supplies auxiliary holomorphic quantities
$h,g$ and two polynomial norm formulas. Applying an exterior-product
identity to the augmented vector $(F_1,h)$ yields
\[
|\langle g,f_0\rangle|^2
=t^{3m-6}\bigl(\mu_2|f_0|^2-c t^m\bigr),\qquad
t=1+z\bar z,\quad m=n-1,
\]
where $c,\mu_2>0$. After complexification, the multiplicity of $1+zw$
forces the inner product to vanish for $m\geq7$. When $m=6$, the
multiplicity bound is sharp, but a polynomial divisibility argument
again gives the required vanishing. Thus the lowest-dimensional case
$n=7$, with $r=1$, is included rather than excluded by an extra
isotropy-order hypothesis.

Related rigidity results concern different additional structures.
Chi, Xie and Xu
\citep{chi2021structure} described two-spheres minimal in both a
hyperquadric and its ambient projective space; simultaneous ambient
minimality is not assumed here. Pacheco and Rehman
\citep[Theorem~3.1 and Corollary~3.2]{pacheco2026primitive} proved
rigidity for full primitive immersions into a full flag manifold having
constant curvature for an invariant flag-manifold metric. Our hypothesis
instead concerns the projected real two-plane map, and does not posit
such a constant-curvature primitive lift. Polynomial norm factorization
also occurs in their argument; here the decisive exterior-product
identity couples the two augmented lifts $(F_1,h)$ and $(F_2,g)$.

Section~\ref{sec:preliminaries} records the geometric conventions and the
construction results used in the proof. Section~\ref{sec:proof} isolates
the polynomial norm input from \citep{li2023conformal}, proves the
exterior-product and divisibility identities, and completes the
classification. Section~\ref{sec:q5} gives the $Q_5$ consequence and a
fully explicit model.

\section{Geometric preliminaries}\label{sec:preliminaries}

\subsection{The quadric and the metric}

The Hermitian product on $\mathbb C^n$ is linear in the first variable:
\[
\langle Z,Y\rangle=\sum_{a=1}^n Z_a\overline{Y_a},\qquad
|Z|^2=\langle Z,Z\rangle.
\]
The complex quadric is
\[
Q_{n-2}=\{[Z]\in\mathbb{CP}^{n-1}:Z^TZ=0\}.
\]
The map
\[
\pi:Q_{n-2}\longrightarrow G(2,n;\mathbb R),\qquad
[Z]\longmapsto\operatorname{span}_{\mathbb R}
\{\operatorname{Re}Z,\operatorname{Im}Z\}
\]
is the double cover which forgets the orientation. We identify a real
plane with its complexification. Under this identification,
$\pi([Z])=\underline Z\oplus\overline{\underline Z}$.
Every sphere in the real Grassmannian lifts to the quadric, since
$\mathbb S^2$ is simply connected.

Represent $G(k,n;\mathbb C)$ by rank-$k$ Hermitian projections $P$ and
use the metric $\frac12\operatorname{tr}(dP\,dP)$. Its restriction
defines the metric on $G(2,n;\mathbb R)$. If $\sigma$ is the inclusion
of the quadric in projective space and $g_{\mathrm{FS}}$ has constant
holomorphic sectional curvature $4$, then
\begin{equation}\label{eq:quadric-metric-normalization}
\pi^*g_{G(2,n;\mathbb R)}=2\sigma^*g_{\mathrm{FS}}.
\end{equation}
This is the normalization used in
\citep{jiao2019q5,li2023conformal}. With the Fubini--Study induced
quadric metric, all the Gaussian curvatures stated here must therefore
be multiplied by two. Multiplication of a target metric by a positive
constant does not change harmonicity or minimality.

Throughout, classification is understood up to target isometries and
conformal reparametrizations of the domain. This convention is needed
when choosing a stereographic coordinate in which a constant-curvature
metric has its standard round form.

\subsection{Harmonic sequences and ramification}

A map $\psi:M\to G(k,n;\mathbb C)$ determines a rank-$k$ subbundle
$\underline\psi\subset M\times\mathbb C^n$. It is holomorphic precisely
when this is a holomorphic subbundle of the trivial holomorphic bundle.
For a local coordinate $z$, write $\partial=\partial/\partial z$ and
$\bar\partial=\partial/\partial\bar z$; primes denote $z$-derivatives
of holomorphic quantities. For a section $s$ of $\underline\psi$, put
\[
A'_\psi s=\operatorname{proj}_{\psi^\perp}(\partial s),\qquad
A''_\psi s=\operatorname{proj}_{\psi^\perp}(\bar\partial s).
\]
For harmonic $\psi$, the Gauss transforms $G'(\psi)$ and $G''(\psi)$
are the saturated images of these forms. Thus their images, initially
defined away from isolated rank drops, are extended across those points
as subbundles. Iteration defines the forward and backward harmonic
sequences $\psi_0=\psi$, $\psi_{j+1}=G'(\psi_j)$ and
$\psi_{-j-1}=G''(\psi_{-j})$; see
\citep{burstall1986construction,bahy1989explicit}.

\begin{definition}\label{def:isotropy}
The isotropy order of $\psi$ is the largest integer $r$ such that
$\underline\psi_j\perp\underline\psi_0$ for every $1\leq j\leq r$.
If all positive transforms are orthogonal to $\underline\psi_0$, set
$r=\infty$. Thus finite isotropy order means that the first failure of
orthogonality occurs at $j=r+1$.
A harmonic map is linearly full if its subbundle is not contained in
any proper constant complex subspace of $\mathbb C^n$. It is
$\partial'$-irreducible if $\operatorname{rank}G'(\psi)
=\operatorname{rank}\psi$. For maps into $G(2,n;\mathbb R)$,
$\partial'$- and $\partial''$-irreducibility agree, and we simply say
irreducible.
\end{definition}

Choose unitary frames $E_j$ of $\underline\psi_j$, and let $\Omega_j$
be the coefficient matrix of $A'_{\psi_j}$, so that
\[
\partial E_j=E_{j+1}\Omega_j+E_j\Psi_j,\qquad
\Omega_j\in M_{k_{j+1}\times k_j}(\mathbb C).
\]
Only arrows with nonzero target rank are considered. Define
$L_j=\operatorname{tr}(\Omega_j\Omega_j^*)$. For the conformal maps
used here, the induced metric and curvature are
\begin{equation}\label{eq:grassmannian-metric}
ds_j^2=(L_{j-1}+L_j)|dz|^2,\qquad
K_j=-\frac{2}{L_{j-1}+L_j}
\partial\bar\partial\log(L_{j-1}+L_j).
\end{equation}
For real $\psi$, the backward sequence is the conjugate of the forward
one; in particular $L_{-1}=L_0$.

\begin{definition}\label{def:unramified}
A forward arrow from $\psi_j$ to its nonzero Gauss transform is
unramified if
\[
\det(\Omega_j\Omega_j^*)|dz|^{2k_{j+1}}
\]
is nowhere vanishing. The analogous condition defines unramified
backward arrows. A harmonic map is totally unramified if every nonzero
arrow in both its forward and backward sequences is unramified,
including the arrows starting at $\psi_0$.
\end{definition}

Here $|dz|^2=dz\,d\bar z$ is a tensor, not a wedge product. The
displayed density is invariant under unitary frame changes and
coordinate changes. The condition requires full row rank of each
$\Omega_j$ everywhere, and is stronger than immersion of $\psi_0$.
For real maps, complex conjugation identifies the two directions, so
it suffices to verify the forward arrows.

\subsection{Projective Frenet frames}\label{sec:projective}

For a linearly full holomorphic curve in $\mathbb{CP}^m$, choose a
polynomial lift $f_0:\mathbb C\to\mathbb C^{m+1}$ whose components have
no common zero. Its orthogonal Frenet vectors are defined away from
ramification points by
\begin{equation}\label{eq:frenet}
f_{j+1}=\partial f_j-(\partial\log|f_j|^2)f_j,\qquad
\bar\partial f_j=-\ell_{j-1}f_{j-1},\qquad
\ell_j=\frac{|f_{j+1}|^2}{|f_j|^2},
\end{equation}
with $\ell_{-1}=\ell_m=0$. The associated osculating lift is
\[
F_j=f_0\wedge f_0'\wedge\cdots\wedge f_0^{(j)}
    =f_0\wedge\cdots\wedge f_j.
\]
If the whole Frenet sequence is totally unramified, these exterior
products have no common zero in their components, and their degrees are
\begin{equation}\label{eq:osculating-degrees}
\deg F_j=(j+1)(m-j),\qquad 0\leq j\leq m.
\end{equation}
In particular $\deg f_0=m$. These are the unramified Pl\"ucker
formulas \citep{bolton1988conformal,li2023conformal}; the hypothesis
concerns the entire sequence, not merely its initial segment.

For reference, start the standard Veronese sequence with
\begin{equation}\label{eq:veronese-seed}
V_0^{(m)}(z)=\left(\sqrt{\binom m0},
\sqrt{\binom m1}z,\ldots,\sqrt{\binom mm}z^m\right)^T
\end{equation}
and define $V_j^{(m)}$ by \eqref{eq:frenet}. For $t=1+|z|^2$,
\begin{equation}\label{eq:veronese-norms}
|V_j^{(m)}|^2=\frac{m!\,j!}{(m-j)!}\,t^{m-2j},\qquad
\ell_j=\frac{(j+1)(m-j)}{t^2}.
\end{equation}
These relations hold in all dimensions considered below.

\subsection{Finite-isotropy construction}

We use the finite-isotropy factorization of Bahy-El-Dien and Wood
\citep[Theorem~4.7]{bahy1989explicit}. For an irreducible real
two-plane map of finite isotropy order, this order is odd, say
$r=2s+1$ with $s\geq0$.
There is an integer $l\geq1$ and a holomorphic initial line
$[f_0]:\mathbb S^2\to\mathbb{CP}^m\subset\mathbb{CP}^{n-1}$,
linearly full in its span, from which the map is obtained by $l$
forward--backward replacement pairs. The intermediate maps take values
in the complex Grassmannian; the even replacement stages are real.
The initial curve satisfies
\begin{equation}\label{eq:bilinear-isotropy}
f_j^Tf_0=0\quad(0\leq j\leq2(l+s)+1),\qquad
f_{2(l+s)+2}^Tf_0\not\equiv0.
\end{equation}
Consequently,
\begin{equation}\label{eq:dimension-bound}
r+2l+1=2(l+s)+2\leq m\leq n-1.
\end{equation}
For $r=n-6$ this implies $l\leq2$. Li and Jiao
\citep[Proposition~3]{li2023conformal} proved that the $l=1$ case
cannot have both total unramification and constant curvature when
$n\geq7$. Their nonexistence result does not assume polynomiality.
Thus every map in Theorem~\ref{thm:main} has
\begin{equation}\label{eq:dimension-equality}
l=2,\qquad m=n-1,\qquad n\ \hbox{odd},\qquad m\geq6.
\end{equation}
In particular, $r=1$ corresponds exactly to $n=7$.

\section{Polynomial rigidity and the proof of the theorem}\label{sec:proof}

\subsection{The two-step norm input}

We recall the two-step construction data of Li and Jiao
\citep[Proposition~2 and equations~(68)--(85)]{li2023conformal}, using
$m=n-1$.

\begin{lemma}[Li--Jiao]\label{lem:setup}
Let $\phi$ satisfy the hypotheses of Theorem~\ref{thm:main}. After a
conformal reparametrization, its initial holomorphic curve admits a
nowhere-zero polynomial lift $f_0:\mathbb C\to\mathbb C^n$ of degree
$m=n-1$, and its entire Frenet sequence is totally unramified. Put
$F_1=f_0\wedge f_0'$, $F_2=f_0\wedge f_0'\wedge f_0''$.
On the open set where the replacement coordinates are defined, write
\begin{align}
v&=f_1+x_0\bar f_0,& h&=x_0|f_0|^2,\label{eq:first-replacement}\\
X&=-\bar x_0|f_0|^2f_1+|f_1|^2\bar f_0,&
C&=\partial x_0-x_0\partial\log|f_1|^2,\nonumber\\
v_1&=f_2+\frac{C}{|v|^2}X+x_1\bar v,&
\phi&=\underline{v_1}\oplus\overline{\underline{v_1}}.
\label{eq:second-replacement}
\end{align}
Then $h$ and the vector
\begin{equation}\label{eq:g-definition}
g=h f_2-\bigl(h'-h\partial\log|F_1|^2\bigr)f_1
       +x_1|v|^2f_0
\end{equation}
extend to polynomials on $\mathbb C$. There are constants
$\lambda_1,\mu_2>0$ such that
\begin{align}
|F_1|^2+|h|^2&=\lambda_1t^{2m-2},\label{eq:first-norm}\\
|F_2|^2+|g|^2&=\mu_2t^{3m-6},\label{eq:second-norm}
\end{align}
where $t=1+|z|^2$. Moreover,
\begin{equation}\label{eq:setup-curvature}
K=\frac{2}{5m-8},\qquad ds^2=\frac{2(5m-8)}{t^2}|dz|^2.
\end{equation}
\end{lemma}

\begin{remark}
The conclusions of Lemma~\ref{lem:setup} are established in
\citep{li2023conformal} before the polynomiality assumption on $|f_2|^2$
is used after equation~(92); that assumption is not needed here.
\end{remark}

\subsection{Exterior-product identities}

All exterior spaces carry the Hermitian metric induced from
$\mathbb C^n$. It is important to distinguish an exterior product in
$\Lambda^2(\Lambda^2\mathbb C^n)$ from multiplication into
$\Lambda^4\mathbb C^n$.

\begin{lemma}\label{lem:exterior}
For the quantities in Lemma~\ref{lem:setup},
\begin{align}
|F_1\wedge F_1'|^2_{\Lambda^2(\Lambda^2\mathbb C^n)}
 &=|f_0|^2|F_2|^2,\label{eq:compound-identity}\\
hF_1'-h'F_1&=f_0\wedge g.\label{eq:auxiliary-exterior}
\end{align}
If a holomorphic polynomial vector $a$ in any Hermitian space satisfies
$|a|^2=\lambda t^N$, then
\begin{equation}\label{eq:round-exterior}
|a\wedge a'|^2=\lambda^2N t^{2N-2}.
\end{equation}
\end{lemma}

\begin{proof}
For \eqref{eq:compound-identity}, apply the Gram determinant identity
to $f_0,f_0',f_0''$. The two-vector Gram determinant of
$f_0\wedge f_0'$ and $f_0\wedge f_0''$ is $|f_0|^2$ times their
three-vector Gram determinant, namely $|f_0|^2|F_2|^2$.

Put $a_0=\partial\log|f_0|^2$ and $a_1=\partial\log|f_1|^2$.
The Frenet equations give
\[
f_1=f_0'-a_0f_0,\qquad
f_2=f_0''-(a_0+a_1)f_0'+(a_0a_1-a_0')f_0.
\]
Since $\partial\log|F_1|^2=a_0+a_1$, substituting into
\eqref{eq:g-definition} yields the exact expansion
\[
g=hf_0''-h'f_0'+qf_0,\qquad
q=a_0h'-h(a_0'+a_0^2)+x_1|v|^2.
\]
Taking the exterior product with $f_0$ proves
\eqref{eq:auxiliary-exterior}.

Finally, for any holomorphic vector $a$, its Gram determinant is
\[
|a\wedge a'|^2
=|a|^2\partial\bar\partial|a|^2
 -(\partial|a|^2)(\bar\partial|a|^2).
\]
Substituting $|a|^2=\lambda t^N$ proves
\eqref{eq:round-exterior}.
\end{proof}

\begin{lemma}\label{lem:fundamental}
Set $N=2m-2$, $D=3m-6$ and $c=\lambda_1^2N$. Then
\begin{equation}\label{eq:fundamental-real}
|\langle g,f_0\rangle|^2=t^D(\mu_2|f_0|^2-ct^m).
\end{equation}
\end{lemma}

\begin{proof}
Let $a=(F_1,h)\in\Lambda^2\mathbb C^n\oplus\mathbb C$, and let $e$
be the unit vector in the last summand. The orthogonal decomposition
\[
a\wedge a'=F_1\wedge F_1'+(h'F_1-hF_1')\wedge e
\]
together with Lemma~\ref{lem:exterior} gives
\begin{align*}
ct^{2N-2}
 &=|f_0|^2|F_2|^2+|f_0\wedge g|^2\\
 &=|f_0|^2\bigl(|F_2|^2+|g|^2\bigr)-|\langle g,f_0\rangle|^2\\
 &=\mu_2|f_0|^2t^D-|\langle g,f_0\rangle|^2.
\end{align*}
Since $2N-2=D+m$, this is \eqref{eq:fundamental-real}.
\end{proof}

\subsection{Multiplicity and the boundary case}

We now work in $\mathbb C[z,w]$, with $w$ independent, and put
$t=1+zw$. For a polynomial define
\[
\left(\sum c_{ab}z^aw^b\right)^\#
   =\sum\overline{c_{ab}}z^bw^a.
\]
On $w=\bar z$ this is complex conjugation. Polynomial identities on
that real slice extend uniquely to $\mathbb C[z,w]$. We use the same
notation for a squared norm and its complexification. The polynomial
$t$ is irreducible, and its multiplicity $\nu_t$ satisfies
$\nu_t(A^\#)=\nu_t(A)$ for every nonzero polynomial $A$.

\begin{lemma}\label{lem:polynomial-divisibility}
Let $m\geq2$, $\alpha,\beta>0$ and $p\in\mathbb C[z]$. If a polynomial
$Q\in\mathbb C[z,w]$ has the form
\[
Q=\alpha t^m+\beta p(z)p^\#(w),
\]
then
\begin{equation}\label{eq:divisibility}
t^{m-2}\mid QQ_{zw}-Q_zQ_w.
\end{equation}
\end{lemma}

\begin{proof}
Put $B=\beta p(z)p^\#(w)$, so that $Q=\alpha t^m+B$.
Then $BB_{zw}-B_zB_w=0$, and direct expansion gives
\begin{align*}
QQ_{zw}-Q_zQ_w
={}&\alpha^2m t^{2m-2}\\
 &+\alpha t^{m-2}\bigl(t^2B_{zw}
       +m(1+mzw)B-mt(wB_w+zB_z)\bigr).
\end{align*}
Both terms are divisible by $t^{m-2}$.
\end{proof}

\begin{proposition}\label{prop:polynomial-rigidity}
Under the hypotheses and notation of Lemma~\ref{lem:setup}, there is
$a>0$ such that
\begin{equation}\label{eq:rigidity-conclusion}
|f_0|^2=at^m,\qquad h=0,\qquad g=0,\qquad x_0=x_1=0.
\end{equation}
The assertions about $x_0,x_1$ hold on the replacement coordinate
domain. In particular, the final replacement line is $\underline{f_2}$.
\end{proposition}

\begin{proof}
Let $A$ be the complexification of $\langle g,f_0\rangle$. Because
$g$ is polynomial and $\deg f_0=m$, we have $\deg_w A\leq m$.
The complexified identity \eqref{eq:fundamental-real} is
\begin{equation}\label{eq:fundamental-polynomial}
AA^\#=t^D(\mu_2|f_0|^2-ct^m),\qquad D=3m-6.
\end{equation}

\emph{The case $m\geq7$.}
If $A\neq0$, then
\[
\nu_t(AA^\#)=2\nu_t(A)\leq2\deg_w A\leq2m<D,
\]
contrary to \eqref{eq:fundamental-polynomial}. Hence $A=0$ and
$|f_0|^2=(c/\mu_2)t^m$.

\emph{The case $m=6$.}
If $A=0$ the same conclusion follows. Suppose instead that $A\neq0$.
Equation~\eqref{eq:fundamental-polynomial} and $\deg_w A\leq6$ give
$12\leq2\nu_t(A)\leq12$. Thus $\nu_t(A)=6$, and division by $t^6$
leaves a polynomial independent of $w$:
\[
A=t^6p(z),\qquad
|f_0|^2=\frac{c}{\mu_2}t^6+\frac{1}{\mu_2}p(z)p^\#(w).
\]
Applying Lemma~\ref{lem:polynomial-divisibility} to $|f_0|^2$ in the
Gram identity
\[
|F_1|^2=|f_0|^2(|f_0|^2)_{zw}-(|f_0|^2)_z(|f_0|^2)_w
\]
gives $t^4\mid|F_1|^2$.
Equation~\eqref{eq:first-norm} then gives $t^4\mid hh^\#$.
A nonzero polynomial in $z$ alone, or in $w$ alone, cannot be divisible
by $1+zw$. Hence $h=0$, and $|F_1|^2=\lambda_1t^{10}$.
Equations~\eqref{eq:compound-identity} and \eqref{eq:round-exterior}
now imply
\begin{equation}\label{eq:boundary-product}
|f_0|^2|F_2|^2=10\lambda_1^2t^{18}.
\end{equation}
Unique factorization in $\mathbb C[z,w]$ forces every irreducible
factor of $|f_0|^2$ to be associated with $t$. The bidegree of $|f_0|^2$
is exactly $(6,6)$, so $|f_0|^2=at^6$ for $a>0$. Substitution into
\eqref{eq:fundamental-polynomial} gives
$AA^\#=(\mu_2a-c)t^{18}$. If the coefficient is nonzero, the right
side has multiplicity $18$, whereas the left side has multiplicity
$12$; if the coefficient is zero, then $A=0$. Both contradict the
assumption $A\neq0$. Thus $A=0$ also in the boundary case.

We have proved in both cases that $|f_0|^2=at^m$ and $A=0$. Applying the
Gram identity to $f_0$ yields $|F_1|^2=a^2m t^{2m-2}$, so
\eqref{eq:first-norm} reads
\[
hh^\#=(\lambda_1-a^2m)t^{2m-2}.
\]
The same one-variable factor argument gives $h=0$. Since $|f_0|^2>0$ on
$\mathbb C$, equation~\eqref{eq:first-replacement} implies $x_0=0$.
Equation~\eqref{eq:g-definition} reduces to
$g=x_1|f_1|^2f_0$, and therefore
$0=\langle g,f_0\rangle=x_1|f_1|^2|f_0|^2$. Total unramification gives
$|f_1|^2>0$, whence $x_1=0$ and $g=0$. Substitution in
\eqref{eq:second-replacement} gives $v_1=f_2$, as required.
The equality of the corresponding line subbundles extends across the
exceptional coordinate points by smoothness.
\end{proof}

\subsection{Completion of the classification}

\begin{proof}[Proof of Theorem~\ref{thm:main}]
By Proposition~\ref{prop:polynomial-rigidity}, $|f_0|^2=at^m$ for
some $a>0$, where $m=n-1\geq6$. The Frenet equations
\eqref{eq:frenet} give
\[
|f_2|^2=2am(m-1)t^{m-4},
\]
which is polynomial in $z,\bar z$. Thus the additional hypothesis of
\citep[Theorem~1]{li2023conformal} is automatic, and that theorem gives
\eqref{eq:main-model}. The converse and the target-isometry uniqueness
are the known properties of these Veronese models; see
\citep[Proposition~4 and Remark~1]{li2023conformal}.
\end{proof}

\section{The boundary model in $Q_5$}\label{sec:q5}

\begin{corollary}\label{cor:q5-without-polynomiality}
A linearly full, irreducible, totally unramified conformal minimal
immersion $\phi:\mathbb S^2\to G(2,7;\mathbb R)$ of constant Gaussian
curvature and isotropy order one is, after a conformal reparametrization
and a target isometry, the model in Example~\ref{ex:q5-model}. Its
curvature is $1/11$. No separate two-step or polynomiality assumption
is needed.
\end{corollary}

\begin{proof}
Apply Theorem~\ref{thm:main} with $n=7$. The uniqueness statement
identifies the model explicitly given below. The elimination of the
polynomiality assumption strengthens the irreducible two-step branch
of \citep[Proposition~4]{jiao2019q5}; it does not subsume the reducible
cases studied in that paper.
\end{proof}

\begin{example}\label{ex:q5-model}
We express the model of \citep[Proposition~4]{jiao2019q5} in the
notation of the proof. Put $t=1+|z|^2$ and
\[
V_0^{(6)}=(1,\sqrt6z,\sqrt{15}z^2,2\sqrt5z^3,
             \sqrt{15}z^4,\sqrt6z^5,z^6)^T,
\]
\[
U=\frac1{\sqrt2}
\begin{pmatrix}
1&0&0&0&0&0&1\\
i&0&0&0&0&0&-i\\
0&1&0&0&0&-1&0\\
0&i&0&0&0&i&0\\
0&0&1&0&1&0&0\\
0&0&i&0&-i&0&0\\
0&0&0&i\sqrt2&0&0&0
\end{pmatrix}.
\]
Then $U^*U=I$ and $U^TU=W_7$. Set $f_0=UV_0^{(6)}$ and
\[
f_{j+1}=\partial f_j-\frac{(6-2j)\bar z}{t}f_j,
\qquad 0\leq j\leq5.
\]
In particular,
\begin{equation}\label{eq:q5-explicit-lift}
f_2=f_0''-\frac{10\bar z}{t}f_0'
                  +\frac{30\bar z^2}{t^2}f_0.
\end{equation}
The polynomial representative in equation~(4.30) of
\citep{jiao2019q5} is $P=\tfrac12t^2f_2$. Thus $P$ and $f_2$
define the same projective map, although their norms differ.
Directly from \eqref{eq:veronese-norms},
\[
|f_j|^2=\frac{6!\,j!}{(6-j)!}t^{6-2j},\qquad 0\leq j\leq6,
\]
and $P^TP=0$, $|P|^2=15t^6$. Moreover,
$\overline{\underline{f_j}}=\underline{f_{6-j}}$, so
$[f_2]\in Q_5$ and $\phi=\underline{f_2}\oplus\underline{f_4}$.
The coefficient vectors of $P$ as a polynomial in $z,\bar z$ span
$\mathbb C^7$. At infinity, with $\zeta=1/z$, the lift
$\widetilde f_2(\zeta)=\zeta^2 f_2(1/\zeta)$ extends smoothly and has
squared norm $60(1+|\zeta|^2)^2$.

The forward harmonic sequence is
\[
\begin{gathered}
\underline{f_2}\oplus\underline{f_4}
\longrightarrow\underline{f_3}\oplus\underline{f_5}
\longrightarrow\underline{f_4}\oplus\underline{f_6}\\
\longrightarrow\underline{f_5}
\longrightarrow\underline{f_6}\longrightarrow0.
\end{gathered}
\]
The first transform is orthogonal to $\phi$ and has rank two; the
second meets $\phi$ in $\underline{f_4}$. Thus the map is irreducible
and its isotropy order, in the consecutive sense of
Definition~\ref{def:isotropy}, is one. The Gram determinants of the
four nonzero arrows are
\[
\det(\Omega_0\Omega_0^*)=\frac{120}{t^4},\qquad
\det(\Omega_1\Omega_1^*)=\frac{72}{t^4},
\]
\[
\det(\Omega_2\Omega_2^*)=\frac{10}{t^2},\qquad
\det(\Omega_3\Omega_3^*)=\frac6{t^2}.
\]
The associated densities have no zeros, including at infinity.
Complex conjugation gives the same conclusion for the backward
sequence. The construction has $x_0=x_1=0$, hence $h=g=0$, and
\[
ds^2=2\left(\frac{|f_2|^2}{|f_1|^2}
                 +\frac{|f_3|^2}{|f_2|^2}\right)|dz|^2
     =\frac{44}{t^2}|dz|^2,\qquad K=\frac1{11}.
\]
For the Fubini--Study induced quadric metric the curvature is $2/11$,
as prescribed by \eqref{eq:quadric-metric-normalization}.

Finally, the exterior products are
\[
|F_1|^2=6t^{10},\qquad |F_2|^2=360t^{12},\qquad
|F_1\wedge F_1'|^2=360t^{18}=|f_0|^2|F_2|^2.
\]
Thus $N=10$, $D=12$, $\lambda_1=6$ and
$\mu_2=360=\lambda_1^2N$. These formulas exhibit the equality of the
multiplicity bounds in the boundary dimension; its rigidity is proved
by Proposition~\ref{prop:polynomial-rigidity}, not by assuming the
polynomiality visible in this already classified example.
\end{example}

\section{Concluding remarks}

Theorem~\ref{thm:main} removes the extra Frenet-norm polynomiality
condition at isotropy order $n-6$. The geometric reduction, the
one-step nonexistence result and the two-step norm input come from the
earlier construction theory and \citep{li2023conformal}; the rigidity
step here is the augmented exterior-product identity and its
polynomial-factor consequence, including the boundary case $n=7$.
The explicit $Q_5$ model makes this boundary mechanism concrete.

The argument retains total unramification. Allowing ramification would
introduce additional factors into the osculating products and change
the degree bounds used in the multiplicity comparison. Similarly,
smaller finite isotropy orders permit further replacement steps.
Neither extension follows from the present proof.

\section*{Acknowledgments}

This work was supported by the National Natural Science Foundation of
China (Grant No.~12371055).

\bibliographystyle{plain}
\bibliography{01}

\section*{AI Assistance Statement}

AI assisted with the writing and revision of this paper.

\end{document}